\documentclass[12pt]{article}
\usepackage{epsfig}
\usepackage{graphicx,psfrag}

\usepackage{multirow, color}

\usepackage{amssymb}
\begin{document}

\begin{titlepage}
\title{Specification Test based on Convolution-type Distribution Function Estimates for Non-linear Auto-regressive Processes}
\author{Kun Ho Kim, and Jiwoong Kim}
\maketitle
\abstract{The paper proposes a specification test based on two estimates of distribution function. One is the traditional kernel distribution function estimate and the other is a newly proposed convolution-type distribution function estimate. Asymptotic properties of the new estimate are studied when the innovation density is known and when it is unknown. The MISE-type statistic based on these estimates is suggested to test parametric specifications of the mean and volatility functions. The relating asymptotic results are obtained and the finite-sample properties are studied based on the bootstrap methodology. A simulation study shows that the proposed test competes favorably to benchmark tests in terms of the empirical level and power.\\}

\vspace{0.2cm}
\end{titlepage}

\section{Introduction}

    Consider the following model framework in time series:
\begin{eqnarray}\label{model0}		
X_i=\mu(X_{i-1})+\sigma(X_{i-1})\epsilon_i,~~~~~~~~i\in\mathbb{Z}
\end{eqnarray}
    where $X_i$ is a stationary process, $\mu:\mathbb R\rightarrow\mathbb{R}$ and $\sigma:\mathbb{R}\rightarrow\mathbb{R}$ are unknown conditional mean and variance functions respectively, and $\epsilon_i$ are  independent and identically distributed (iid) innovations. The purpose of this paper is to construct a specification test based on cumulative distribution function (c.d.f.) estimates for $X_i$ in (\ref{model0}).
		
		Various tests for (\ref{model0}) have been proposed in the time series literature: A\"{i}t-Sahalia (1996) proposed a parametric specification test by comparing the nonparametric kernel density estimate of the marginal density of $X_i$ with its closed-form density estimate under the parametric form. Given that the kernel density estimate always converges to the true density, the difference between these estimates would converge to zero only if the parametric forms are correctly specified. In stead of using the density estimates, Corradi and Swanson (2005) propose a test that utilizes the empirical c.d.f. of $X_i$ and the closed-form c.d.f. estimate under the parametric forms of mean and variance functions. Given that the limiting distribution of their test statistic is a functional of a Gaussian process, they employ bootstrap procedures to carry out inference.
		
	  Despite the innovative nature of the idea, the tests in A\"{i}t-Sahalia (1996) and Corradi and Swanson (2005) are not applicable if the closed-form density and closed-form c.d.f. of $X_i$ are not available. This significantly reduces the applicability of these tests because the closed-form density and c.d.f. are typically unavailable for many prominent non-linear time series models, such as the autoregressive conditional heteroskedastic (ARCH) process. To address this issue, Kim, Zhang and Wu (2015) introduce a convolution-type {\it density} estimate that is used to test for the framework (\ref{model0}). They propose a test statistic based on the maximal deviation of this convolution density estimate from the traditional kernel density estimate. Since the convolution only requires independence between the mean and variance, one needs not to know the closed-form density/c.d.f. of $X_i$, which greatly enhances the applicability of the proposed test.
		
		The potential problem of the test in Kim et al (2005) is that their test depends on the kernel density estimate of which the convergence is very slow. The slow convergence could potentially lead to size distortion and a low power of the test. One way to address this issue is to construct a test that employs c.d.f. estimates based on kernel smoothing and convolution. Given the additive mean and multiplicative variance of (\ref{model0}), the convolution can be applied to obtain a c.d.f. estimate of $X_i$. A simple modification of (\ref{model0}) shows:
\begin{eqnarray}\label{model01}		
Y_i=m(X_{i-1})+\epsilon_i
\end{eqnarray}
    where $Y_i:=X_i/\sigma(X_{i-1})$ and $m(x):=\mu(x)/\sigma(x)$. Note that there is independence between $m(X_{i-1})$ and $\epsilon_i$, such that the convolution applies. While Kim et al (2015) use it for {\it density} estimation, we use convolution to estimate the c.d.f. of $X_i$ and compare it to the kernel c.d.f estimate because both the convolution-type c.d.f. estimate and the kernel c.d.f. estimate achieve the root-n-consistency. As shown in Remark 1 of this paper, our test based on the kernel c.d.f. estimate and the convolution-type c.d.f. estimate enjoys a faster convergence than that in Kim et al (2015), which explains why our test tends to perform better than its competitors including that in Kim et al (2015), as shown by the simulation study in this work.
		
		The organization of the paper is the following: Section 2 introduces the technical assumptions required for our results and discusses the methodology on the convolution-based c.d.f. estimate. We first consider the case of known innovation density to derive the relating asymptotic properties. The result is later generalized to the case of unknown innovation density, which has more of practical relevance than the former. Section 3 describes how to construct a specification test based on the kernel c.d.f. estimate and the convolution c.d.f. estimate from the previous section. The asymptotic distribution of the test statistics is derived, and its finite-sample properties based on the {\it bootstrap} methodology are examined by a simulation study. Section 4 concludes the paper and discusses related future research. Tables and figures are relegated to the appendix of the paper.

\section{Methodology}

	  For simplicity, we first consider the autoregressive process with a homoskedastic innovation:
\begin{eqnarray}\label{model1}		
X_i=m_{\theta}(X_{i-1})+\epsilon_i
\end{eqnarray}
    where $\theta$ is an unknown parameter. The homoskedasticity assumption here will be relaxed to the case of conditional heteroskedasticity later. The parametric specification that needs to be tested is the following:
\begin{eqnarray}\label{null}		
H_0: m(\cdot)=m_{\theta}(\cdot)
\end{eqnarray}	
    where $m(\cdot)$ is the true mean function. Testing (\ref{null}) has been conducted in various contexts. Corradi and Swanson (2005) propose a MISE-type statistic that utilizes an empirical c.d.f. of $X_i$ and its parametric closed-form c.d.f. function. Kim et al (2015) consider a kernel density estimate and a convolution-type density estimate for $X_i$ to test (\ref{null}) using the maximal distance between the two estimates.
		
		In this paper, we combine the ideas of Corradi and Swanson (2005) and Kim et al (2015) to propose the test statistic based on the kernel c.d.f. estimate and the convolution-type c.d.f. estimate for $X_i$. First, define the kernel c.d.f. estimate for $X_i$:
\begin{eqnarray}\label{kern}
\hat{F}_{k}(x)=\frac{1}{n}\sum_{i=1}^nG\left(\frac{x-X_i}{b}\right)
\end{eqnarray}
    where $G(u)=\int_{-\infty}^uK(x)dx$ and $K(\cdot)$ is a kernel function. Here $b$ is a bandwidth. Given the independence between the mean and innovation of (\ref{model1}), the c.d.f. of $X_i$, $F_X(\cdot)$, also can be estimated by the following convolution c.d.f. estimate:
\begin{eqnarray}\label{convol}
\hat{F}_{c}(x)=\int_{\mathbb{R}}\hat{F}_{\epsilon}(x-t)\hat{f}_g(t)dt
\end{eqnarray}
   where $F_{\epsilon}(\cdot)$, the c.d.f. of $\epsilon_i$, and $f_g(\cdot)$, the density function of $m_{\theta}(X_{i-1})$, are estimated, respectively, by:
\begin{eqnarray*}
\hat{F}_{\epsilon}(x)=\frac{1}{n}\sum_{i=1}^nG\left(\frac{x-\hat{\epsilon}_i}{b}\right)
\end{eqnarray*}
\begin{eqnarray*}
\hat{f}_g(x)=\frac{1}{nb}\sum_{i=1}^nK\left(\frac{x-m_{\hat{\theta}}\left(X_{i-1}\right)}{b}\right)
\end{eqnarray*}
   Here $\hat{\epsilon}_i=X_i-m_{\hat\theta}(X_{i-1})$ and $\hat{\theta}$ is an $\sqrt{n}$-consistent estimator of parameter $\theta$ in (\ref{model1}), respectively. Obviously, $\hat{F}_{k}(x)$ in (\ref{kern}) converges to the true c.d.f. of $X_i$ regardless of the parametric form of $m_{\theta}(\cdot)$ in (\ref{model1}), while $\hat{F}_{c}(x)$ in (\ref{convol}) converges only under its correct form. Hence the properly centered and scaled difference between the two c.d.f. estimates in (\ref{kern}) and (\ref{convol}) can be used as a statistic for testing the parametric specification of (\ref{model1}).

\subsection{Assumptions}
	
	 Some notations are needed to introduce the assumptions in this study. For a random variable $W$, write $W \in \mathcal{L}^p$, $p > 0$, if $\| W \|_p := [\mathbb{E}(|W|^p)]^{1/p} < \infty$, and write $\| W \| = \| W \|_2$. We define the projection operator $\mathcal{P}$ as $\mathcal{P}_i[\cdot] \equiv \mathbb{E} [\cdot |\mathcal{\cal I}_i] - \mathbb{E} [\cdot |\mathcal{\cal I}_{i-1}]$, where ${\cal I}_i = (\epsilon_i, \epsilon_{i-1},\ldots )$. The following assumptions are needed to derive the asymptotic properties of the convolution c.d.f. estimator:\\
		
{\bf Assumption 1.} {\it Let the kernel function $K$ be bounded, symmetric, with bounded support $[-A, A]$, $K \in {\cal C}^1[-A, A]$, $K(\pm A) = 0$ and $\sup_u |K'(u)| < \infty$.}
		
{\bf Assumption 2.} {\it $\sup_{x \not= x'} |m_{\theta}(x) - m_{\theta}(x')| / |x-x'| < 1$, and $\epsilon_i \in {\cal L}^p$, $p > 0$}

{\bf Assumption 3.} {\it $\sup_{x} \left[ F_{\epsilon}(x) + \left|F_{\epsilon}^{'}(x)\right| + \left|F_{\epsilon}^{''}(x)\right| \right] < \infty$, and as $|x| \to \infty$, $F_{\epsilon}(x) = O(|x|^{-\beta})$ for some $\beta > 0$.}

{\bf Assumption 4.} {\it $\hat{\theta}$ is an estimate of $\theta$ such that
\begin{eqnarray}\label{Bhd}
\sqrt{n}\left(\hat\theta - \theta\right) = \frac{1}{\sqrt{n}}\sum_{i=1}^n R_i + o_\mathbb{P}\left(1\right)
\end{eqnarray}
where $R_i = R(\epsilon_i, \epsilon_{i-1}, \ldots )$ satisfies the short-range dependence condition
\begin{eqnarray}\label{srd}
\sum_{i=0}^\infty \|{\cal P}_0 R_i\| < \infty.
\end{eqnarray}}

{\bf Assumption 5.} {\it Let $\dot{m}_{\theta}(x)=\partial m_{\theta}(x)/\partial\theta$ exist, and
$\left|m_{\theta}(x)-m_{\theta_0}(x)\right|\leq\dot{m}_{\theta}(x)|\theta-\theta_0|$ with $\mathbb{E}\left[\dot{m}_{\theta}^2(X_i)\right]<\infty$.}\\

    Assumption 1 allows popular kernels such as Parzen, Epanechnikov and uniform kernels among others. Assumption 2 represents a contraction condition and it ensures that process $X_i$ is a stationary and ergodic solution of the form $X_i=G\left(\epsilon_i,\epsilon_{i-1},\cdots\right)$. The process is also causal. For many non-linear times series models, the innovation c.d.f. satisfies Assumption 3. Assumption 4 is an important intermediate step in obtaining a central limit theorem for an estimate $\hat\theta$. In certain situations, (\ref{Bhd}) is called the Bahadur representation. Assumption 5 is not the weakest possible. Based on these assumptions, we introduce the convolution c.d.f. estimate and investigate its asymptotic properties.

\subsection{Convolution c.d.f. estimation}

    Let $\hat{S}_n(x)=\sum_{i=1}^nF_{\epsilon}\left(x-m_{\hat\theta}(X_{i-1})\right)$, where $\hat\theta$ is an $\sqrt{n}$-consistent estimate of $\theta$ and $F_{\epsilon}(\cdot)$ is the c.d.f. of $\epsilon_i$ in (\ref{model1}). Given $F_{\epsilon}(\cdot)$, the convolution c.d.f. estimator is:
\begin{eqnarray}\label{kd}
\breve{F}_c(x)&=&\int_{\mathbb{R}}\frac{1}{nb}\sum_{i=1}^nK\left(\frac{t-m_{\hat\theta}(X_{i-1})}{b}\right)F_{\epsilon}(x-t)dt\nonumber\\
              &=&\frac{1}{n}\int\,K(u)\sum_{i=1}^nF_{\epsilon}\left(x-ub-m_{\hat\theta}(X_{i-1})\right)du\nonumber\\
              &=&\frac{1}{n}\int\,K(u)\hat{S}_n(x-ub)du
\end{eqnarray}
    Given Assumptions 1--5, we introduce the following lemmas:\\
					
{\bf Lemma 1.} Given the c.d.f. of $\epsilon_i$,  $F_{\epsilon}(\cdot)$,
\begin{eqnarray}\label{lemma1}
\frac{1}{\sqrt{n}}\left(\hat{S}_n(x)-nF_X(x)\right)\,\,\Rightarrow\,\,N\left(0,\,\,\sigma_1^2(x)\right)
\end{eqnarray}
  where $F_X(\cdot)$ is the c.d.f. of $X_i$, and $\sigma_1(x)=\left\|\sum_{i=1}^{\infty}\mathcal{P}_0\left[F_{\epsilon}\left(x-m_{\theta}(X_{i-1})\right)+c_0R\left(X_{i-1},\theta,\epsilon_i\right)\right]\right\|$.
	
proof) Note that Assumption 4 ensures:
\begin{eqnarray}\label{thest}
\sqrt{n}\left(\hat\theta-\theta\right)=-\frac{1}{\sqrt{n}}\sum_{i=1}^n\,R\left(X_{i-1},\theta,\epsilon_i\right)+o_{\mathbb{P}}(1)
\end{eqnarray}
  where $\mathbb{E}R\left(X_{i-1},\theta,\epsilon_i\right)=0$. Moreover, by the ergodicity of $X_i$ under Assumption 2,
\begin{eqnarray}\label{erg}
c_n=\frac{1}{n}\sum_{i=1}^nf_{\epsilon}\left(x-m_{\theta}(X_{i-1})\right)\,\dot{m}_{\theta}\,\,\,\,\,\,\,\stackrel{\mathbb{P}}{\rightarrow}\,\,\,\,\,\,\,c_0=\mathbb{E}\left[f_{\epsilon}\left(x-m_{\theta}(X_{i-1})\right)\dot{m}_{\theta}\right]
\end{eqnarray}
  where $\dot{m}_{\theta}(x)=\partial m_{\theta}(x)/\partial\theta$ and $f_{\epsilon}(\cdot)$ is the density function for $\epsilon_i$. By a Taylor's expansion of $\hat{S}_n(x)$,
\begin{eqnarray}\label{te1}
\hat{S}_n(x)=\sum_{i=1}^n\,F_{\epsilon}\left(x-m_{\theta}(X_{i-1})\right)-\left(\hat\theta-\theta\right)\sum_{i=1}^nf_{\epsilon}\left(x-m_{\theta}(X_{i-1})\right)\,\dot{m}_{\theta}+O_{\mathbb{P}}\left((\hat\theta-\theta)^2\right)
\end{eqnarray}
  Then, by (\ref{thest})--(\ref{te1}),
\begin{eqnarray}\label{f4}
\frac{1}{\sqrt{n}}\left(\hat{S}_n(x)-nF_X(x)\right)&=&\frac{1}{\sqrt{n}}\sum_{i=1}^n\left[F_{\epsilon}\left(x-m_{\theta}(X_{i-1})\right)-F_X(x)+c_0\,R(X_{i-1},\theta,\epsilon_i)\right]-c_n\,o_{\mathbb{P}}(1)\nonumber\\
&&~~~~+\frac{1}{\sqrt{n}}O_{\mathbb{P}}\left((\hat\theta-\theta)^2\right)+(c_n-c_0)\frac{1}{\sqrt{n}}\sum_{i=1}^n R\left(X_{i-1},\theta,\epsilon_i\right)
\end{eqnarray}
   By Theorem 2 in Wu and Shao (2004),
\begin{eqnarray}\label{f2}
&&\sum_{i=1}^n\left\|\mathcal{P}_0\left[F_{\epsilon}\left(x-m_{\epsilon}(X_{i-1})\right)-F_X(x)+c_0R\left(X_{i-1},\theta,\epsilon_i\right)\right]\right\|\nonumber\\
&\leq&\sum_{i=1}^n\left\|\mathcal{P}_0F_{\epsilon}\left(x-m_{\epsilon}(X_{i-1})\right)\right\|+c_0\sum_{i=1}^n\left\|\mathcal{P}_0R\left(X_{i-1},\theta,\epsilon_i\right)\right\|<\infty
\end{eqnarray}
   Then, by (\ref{f2}), $\mathbb{E}R\left(X_{i-1},\theta,\epsilon_i\right)=0$, $\mathbb{E}\left[F_{\epsilon}\left(x-m_{\theta}(X_{i-1})\right)-F_X(x)\right]=0$ and Theorem 3 in Wu (2005),
\begin{eqnarray}\label{f3}
\frac{1}{\sqrt{n}}\sum_{i=1}^n\left[F_{\epsilon}\left(x-m_{\theta}(X_{i-1})\right)-F_X(x)+c_0\,R(X_{i-1},\theta,\epsilon_i)\right]\,\,\Rightarrow\,\,N\left(0,\,\,\sigma_1^2(x)\right)
\end{eqnarray}
   By applying (\ref{thest}), (\ref{erg}), (\ref{f3}) and  $\frac{1}{\sqrt{n}}\sum_{i=1}^nR\left(X_{i-1},\theta,\epsilon_i\right)\,\rightarrow\,N\left(0,\|\sum_{i=1}^{\infty}\mathcal{P}_0R\left(X_{i-1},\theta,\epsilon_i\right)\|^2\right)$ to (\ref{f4}), the lemma follows.
$$\qquad\qquad\qquad\qquad\qquad\qquad\qquad\qquad\qquad\qquad\qquad\qquad\qquad\qquad\qquad\qquad\qquad\Box$$

   For Lemmas 2--4, we define the following processes:
\begin{eqnarray*}
\hat{F}_{\epsilon}\left(x,\hat\theta\right)&=&\frac{1}{n}\sum_{i=1}^nG\left(\frac{x-{X}_i+m_{\hat\theta}(X_{i-1})}{b}\right)\\
\hat{F}_{\epsilon}\left(x,\theta_0\right)&=&\frac{1}{n}\sum_{i=1}^nG\left(\frac{x-{X}_i+m_{\theta_0}(X_{i-1})}{b}\right)\\
\hat{f}_{g}\left(x,\hat\theta\right)&=&\frac{1}{nb}\sum_{i=1}^nK\left(\frac{x-m_{\hat\theta}\left(X_{i-1}\right)}{b}\right)\\
\hat{f}_{g}\left(x,\theta_0\right)&=&\frac{1}{nb}\sum_{i=1}^nK\left(\frac{x-m_{\theta_0}\left(X_{i-1}\right)}{b}\right)
\end{eqnarray*}
   where $G(u)=\int_{-\infty}^uK(x)dx$.

{\bf Lemma 2.}
\begin{eqnarray*}
\int_{\mathbb{R}}\left(\hat{F}_{\epsilon}\left(x,\hat\theta\right)-\hat{F}_{\epsilon}\left(x,\theta_0\right)\right)^2dx=O_{\mathbb{P}}\left(\frac{1}{nb}\right)
\end{eqnarray*}

proof) Let $C>0$ such that $\int_{\mathbb{R}}|G(u+\delta)-G(u)|^2du\leq C\delta^2$. By the Cauchy-Schwarz inequality,
\begin{eqnarray*}
&&\int_{\mathbb{R}}\left(\hat{F}_{\epsilon}\left(x,\hat\theta\right)-\hat{F}_{\epsilon}\left(x,\theta_0\right)\right)^2dx\\
&=&\frac{1}{n^2}\int_{\mathbb{R}}\left[\sum_{i=1}^n\left(G\left(\frac{x-\hat{\epsilon}_i}{b}\right)-G\left(\frac{x-\epsilon_i}{b}\right)\right)\right]^2dx\\
&\leq&\frac{1}{n^2}\int_{\mathbb{R}}n\sum_{i=1}^n\left[G\left(\frac{x-\hat{\epsilon}_i}{b}\right)-G\left(\frac{x-\epsilon_i}{b}\right)\right]^2dx\\
&=&\frac{1}{n}\sum_{i=1}^n\int_{\mathbb{R}}\left[G\left(\frac{x-\hat{\epsilon}_i}{b}\right)-G\left(\frac{x-\epsilon_i}{b}\right)\right]^2dx\\
&=&\frac{1}{n}\sum_{i=1}^n\int_{\mathbb{R}}\left[G\left(u+\frac{\epsilon_i-\hat\epsilon_i}{b}\right)-G(u)\right]^2bdu\,\leq\,\frac{Cb}{n}\sum_{i=1}^n\left(\frac{\epsilon_i-\hat\epsilon_i}{b}\right)^2
\end{eqnarray*}
   By Assumptions 4 and 5,
\begin{eqnarray*}
\int_{\mathbb{R}}\left(\hat{F}_{\epsilon}\left(x,\hat\theta\right)-\hat{F}_{\epsilon}\left(x,\theta_0\right)\right)^2dx
&\leq&\frac{Cb}{n}\sum_{i=1}^n\left(\frac{m_{\hat\theta}\left(X_{i-1}\right)-m_{\theta_0}\left(X_{i-1}\right)}{b}\right)^2\\
&\leq&\frac{Cb}{n}\sum_{i=1}^n\left(\frac{\left|\hat\theta-\theta_0\right|\dot{m}_{\theta_0}(X_{i-1})}{b}\right)^2=O_{\mathbb{P}}\left(\frac{1}{nb}\right)
\end{eqnarray*}
$$\qquad\qquad\qquad\qquad\qquad\qquad\qquad\qquad\qquad\qquad\qquad\qquad\qquad\qquad\qquad\qquad\qquad\Box$$

{\bf Lemma 2-1.}
\begin{eqnarray}\label{ex2-1}
\int_{\mathbb{R}}\left(\hat{F}_{\epsilon}\left(x,\theta_0\right)-\mathbb{E}\hat{F}_{\epsilon}\left(x,\theta_0\right)\right)^2dx=O_{\mathbb{P}}(b/n)
\end{eqnarray}

proof) Let $l_{i}(x):= G\left((x-\epsilon_i)/b\right)$. Note that
\begin{eqnarray}
\mathbb{E}\left( \hat{F}_{\epsilon}\left(x,\theta_0\right)-\mathbb{E}\hat{F}_{\epsilon}\left(x,\theta_0\right)\right)^2 &=& \mathbb{E} \left\{ \frac{1}{n}\sum_{i=1}^{n}\left( l_{i}(x)-\mathbb{E}l_{i}(x)\right) \right\}^2\nonumber\\
&=& \frac{1}{n^2} \sum_{i=1}^{n} \mathbb{E}(l_{i}(x) - \mathbb{E}l_{i}(x) )^{2}\nonumber\\
&=& \frac{1}{n} \mathbb{E}(l_{0}(x) - \mathbb{E}l_{0}(x) )^{2}\nonumber\\
&\leq & \frac{1}{n} \mathbb{E}l_{0}^{2}(x).\label{eq-21}
\end{eqnarray}
The second inequality follows from the fact that expectation of the cross product terms is zero due to the independence of $\epsilon_i$. Observe that
\begin{eqnarray}
\mathbb{E}l_{0}^{2}(x)&=& \int_{\mathbb{R}} G^{2}\left(\frac{x-y}{b}\right)  dF_{\epsilon}(y)\nonumber\\
&=& b\int_{\mathbb{R}}\int_{\mathbb{R}}G^{2}(u)\,du\, dF_{\epsilon}(y)\nonumber\\
&=& O(b). \label{eq-22}
\end{eqnarray}
By Fubini's theorem, (\ref{eq-21}), and (\ref{eq-22}), (\ref{ex2-1}) follows, thereby completing the proof of the lemma.

{\bf Lemma 3.}
\begin{eqnarray*}
\int_{\mathbb{R}}\left(\hat{F}_{\epsilon}\left(x,\theta_0\right)-F_{\epsilon}(x)\right)^2dx=O_{\mathbb{P}}\left(b^4+\frac{1}{nb}\right)
\end{eqnarray*}

proof) Note the following:
\begin{eqnarray}\label{eq36}
\mathbb{E}\hat{F}_{\epsilon}(x,\theta_0)-F_{\epsilon}(x)&=&\mathbb{E}G\left(\frac{x-\epsilon_i}{b}\right)-F_{\epsilon}(x)\nonumber\\
                                                        &=&b\int_{\mathbb{R}}G(u)f_{\epsilon}(x-ub)du-F_{\epsilon}(x)\nonumber\\
                                                        &=&b\left(-\frac{1}{b}G(u)F_{\epsilon}(x-ub)\Big|^{u=\infty}_{u=-\infty}+\frac{1}{b}\int_{\mathbb{R}}K(u)F_{\epsilon}(x-ub)du\right)-F_{\epsilon}(x)\nonumber\\
                                                        &=&\int_{\mathbb{R}}K(u)F_{\epsilon}(x-ub)du-F_{\epsilon}(x)
\end{eqnarray}
   By a Taylor's expansion on (\ref{eq36}),
\begin{eqnarray}\label{ex1}
\mathbb{E}\hat{F}_{\epsilon}(x,\theta_0)-F_{\epsilon}(x)=\frac{b^2\phi_K}{2}f^{\prime}_{\epsilon}(x)+\frac{b^4}{4!}\int_{\mathbb{R}}u^4K(u)f_{\epsilon}^{(3)}(x_0)du
\end{eqnarray}
   where $x_0\in(x,\,x-ub)$. Moreover,
\begin{eqnarray}\label{ex2}
\int_{\mathbb{R}}\left(\hat{F}_{\epsilon}\left(x,\theta_0\right)-\mathbb{E}\hat{F}_{\epsilon}\left(x,\theta_0\right)\right)^2dx=O_{\mathbb{P}}(1/(nb))
\end{eqnarray}
   By (\ref{ex1}) and (\ref{ex2}),
\begin{eqnarray*}
&&\int_{\mathbb{R}}\left(\hat{F}_{\epsilon}\left(x,\theta_0\right)-F_{\epsilon}(x)\right)^2dx\\
&=&\int_{\mathbb{R}}\left(\hat{F}_{\epsilon}\left(x,\theta_0\right)-\mathbb{E}\hat{F}_{\epsilon}\left(x,\theta_0\right)+\mathbb{E}\hat{F}_{\epsilon}\left(x,\theta_0\right)-F_{\epsilon}(x)\right)^2dx\\
&\leq&2\int_{\mathbb{R}}\left(\hat{F}_{\epsilon}\left(x,\theta_0\right)-\mathbb{E}\hat{F}_{\epsilon}\left(x,\theta_0\right)\right)^2dx+2\int_{\mathbb{R}}\left(\mathbb{E}\hat{F}_{\epsilon}\left(x,\theta_0\right)-F_{\epsilon}(x)\right)^2dx=O_{\mathbb{P}}\left(\frac{1}{nb}+b^4\right)
\end{eqnarray*}
$$\qquad\qquad\qquad\qquad\qquad\qquad\qquad\qquad\qquad\qquad\qquad\qquad\qquad\qquad\qquad\qquad\qquad\Box$$

{\bf Lemma 4.}
\begin{eqnarray*}
\sup_{x}\left|\int_{\mathbb{R}}\left[\hat{F}_{\epsilon}\left(x-y,\hat{\theta}\right)-F_{\epsilon}(x-y)\right]\,\left[\hat{f}_{g}\left(y,\hat{\theta}\right)-f_{g}(y)\right]dy\right|=O_{\mathbb{P}}\left(\sqrt{\frac{1}{n^2b^4}+\frac{b}{n}}\right)
\end{eqnarray*}
proof) Choose $C_0>0$ such that $\int_{\mathbb{R}}|K(u+\delta)-K(u)|^2du\leq\delta^2C_0$. Then, by the Cauchy-Schwarz inequality and by Assumptions 4 and 5,
\begin{eqnarray}\label{1127eq1}
&&\int_{\mathbb{R}}\left[\hat{f}_{g}\left(y,\hat{\theta}\right)-\hat{f}_{g}\left(y,\theta_0\right)\right]^2dy\nonumber\\
&\leq&\frac{1}{nb^2}\sum_{i=1}^n\int_{\mathbb{R}}\left[K\left(\frac{y-m_{\hat\theta}(X_{i-1})}{b}\right)-K\left(\frac{y-m_{\theta_0}(X_{i-1})}{b}\right)\right]^2dy\nonumber\\
&\leq&\frac{1}{nb^2}\sum_{i=1}^nC_0b\left(\frac{m_{\hat\theta}(X_{i-1})-m_{\theta_0}(X_{i-1})}{b}\right)^2\nonumber\\
&\leq&\frac{1}{nb^2}\sum_{i=1}^nC_0b\left(\frac{|\hat\theta-\theta_0|\dot{m}_{\theta_0}(X_{i-1})}{b}\right)^2\nonumber\\
&=&O_{\mathbb{P}}\left(\frac{1}{nb^3}\right)
\end{eqnarray}
    Then,
\small{\begin{eqnarray*}
&&\left(\int_{\mathbb{R}}\left[\hat{F}_{\epsilon}\left(x-y,\hat{\theta}\right)-F_{\epsilon}(x-y)\right]\,\left[\hat{f}_{g}\left(y,\hat{\theta}\right)-f_{g}(y)\right]dy\right)^2\\
&\leq&\int_{\mathbb{R}}\left[\hat{F}_{\epsilon}\left(x-y,\hat{\theta}\right)-F_{\epsilon}(x-y)\right]^2dy\,\int_{\mathbb{R}}\left[\hat{f}_{g}\left(y,\hat{\theta}\right)-f_{g}(y)\right]^2dy\\
&=&\int_{\mathbb{R}}\left[\hat{F}_{\epsilon}\left(x-y,\hat{\theta}\right)-\hat{F}_{\epsilon}\left(x-y,\theta_0\right)+\hat{F}_{\epsilon}\left(x-y,\theta_0\right)-F_{\epsilon}(x-y)\right]^2dy\\
&&~~~~~~~~~~~~~~~~~~\times\int_{\mathbb{R}}\left[\hat{f}_{g}\left(y,\hat{\theta}\right)-\hat{f}_{g}\left(y,\theta_0\right)+\hat{f}_{g}\left(y,\theta_0\right)-f_{g}(y)\right]^2dy\\
&\leq&\left(2\int_{\mathbb{R}}\left[\hat{F}_{\epsilon}\left(x-y,\hat{\theta}\right)-\hat{F}_{\epsilon}\left(x-y,\theta_0\right)\right]^2dy+2\int_{\mathbb{R}}\left[\hat{F}_{\epsilon}\left(x-y,\theta_0\right)-F_{\epsilon}\left(x-y\right)\right]^2dy\right)\\
&&~~~~\times\left(2\int_{\mathbb{R}}\left[\hat{f}_{g}\left(y,\hat{\theta}\right)-\hat{f}_{g}\left(y,\theta_0\right)\right]^2dy+2\int_{\mathbb{R}}\left[\hat{f}_{g}\left(y,\theta_0\right)-f_{g}\left(y\right)\right]^2dy\right)
\end{eqnarray*}}
    Then, by Lemmas 2 and 3, (\ref{1127eq1}), and Lemma 1 in Kim et al (2015),
\begin{eqnarray}\label{eq38}
\left(\int_{\mathbb{R}}\left[\hat{F}_{\epsilon}\left(x-y,\hat{\theta}\right)-F_{\epsilon}(x-y)\right]\,\left[\hat{f}_{g}\left(y,\hat{\theta}\right)-f_{g}(y)\right]dy\right)^2=O_{\mathbb{P}}\left(\frac{1}{n^2b^4}+\frac{b}{n}\right)
\end{eqnarray}
    Hence the lemma follows.
$$\qquad\qquad\qquad\qquad\qquad\qquad\qquad\qquad\qquad\qquad\qquad\qquad\qquad\qquad\qquad\qquad\qquad\Box$$

{\bf Lemma 5.} Recall $c_0=\mathbb{E}\left[\dot{m}_{\theta}\,f_{\epsilon}\left(x-m_{\theta}\left(X_{i-1}\right)\right)\right]$. Then,
\begin{eqnarray*}
\frac{1}{\sqrt{n}}\sum_{i=1}^n\left[F_{\epsilon}\left(x-m_{\theta}\left(X_{i-1}\right)\right)+F_g\left(x-\epsilon_i\right)-2F_X(x)+c_0R\left(X_{i-1},\theta,\epsilon_i\right)\right]\,\,\Rightarrow\,\,N\left(0,\,\sigma_2^2(x)\right)
\end{eqnarray*}
  where $\sigma_2(x)=\left\|\sum_{i=1}^{\infty}\mathcal{P}_0\left[F_{\epsilon}\left(x-m_{\theta}(X_{i-1})\right)+F_g\left(x-\epsilon_i\right)+c_0R\left(X_{i-1},\theta,\epsilon_i\right)\right]\right\|$.

proof) Recall from (\ref{thest}) that $\mathbb{E}R\left(X_{i-1},\theta,\epsilon_i\right)=0$. Then,
\begin{eqnarray*}
\mathbb{E}\left[F_{\epsilon}\left(x-m_{\theta}\left(X_{i-1}\right)\right)+F_g\left(x-\epsilon_i\right)-2F_X(x)+c_0R\left(X_{i-1},\theta,\epsilon_i\right)\right]=0
\end{eqnarray*}
   By a similar argument to (\ref{f2}),
\begin{eqnarray*}
&&\sum_{i=1}^{\infty}\left\|\mathcal{P}_0\left[F_{\epsilon}\left(x-m_{\epsilon}(X_{i-1})\right)+F_g\left(x-\epsilon_i\right)-2F_X(x)+c_0R\left(X_{i-1},\theta,\epsilon_i\right)\right]\right\|<\infty
\end{eqnarray*}
   Hence, by Theorem 3 in Wu (2005), the lemma follows.
$$\qquad\qquad\qquad\qquad\qquad\qquad\qquad\qquad\qquad\qquad\qquad\qquad\qquad\qquad\qquad\qquad\qquad\Box$$

{\bf Lemma 6.} Define $\tilde{S}_n(x)=\sum_{i=1}^n\left[F_{\epsilon}\left(x-m_{\hat\theta}\left(X_{i-1}\right)\right)+F_g\left(x-\hat\epsilon_i\right)\right]$. Then,
\begin{eqnarray*}
\frac{1}{\sqrt{n}}\left(\tilde{S}_n(x)-2nF_X(x)\right)\,\,\Rightarrow\,\,N\left(0,\,\sigma_2^2(x)\right)
\end{eqnarray*}

proof) By a Taylor's expansion,
\begin{eqnarray}\label{new3}
&&\frac{1}{\sqrt{n}}\sum_{i=1}^n\left[F_{\epsilon}\left(x-m_{\hat\theta}\left(X_{i-1}\right)\right)+F_g\left(x-\hat{\epsilon}_i\right)\right]-\frac{1}{\sqrt{n}}\sum_{i=1}^n\left[F_{\epsilon}\left(x-m_{\theta}\left(X_{i-1}\right)\right)+F_g\left(x-\epsilon_i\right)\right]\nonumber\\
&=&\left(\hat\theta-\theta\right)\frac{1}{\sqrt{n}}\frac{\partial}{\partial\theta}\sum_{i=1}^n\left[F_{\epsilon}\left(x-m_{\theta}\left(X_{i-1}\right)\right)+F_g\left(x-\epsilon_i\right)\right]+\frac{1}{\sqrt{n}}O\left(\left(\hat\theta-\theta\right)^2\right)\nonumber\\
&=&-\sqrt{n}\left(\hat\theta-\theta\right)c_n+\frac{1}{\sqrt{n}}O\left(\left(\hat\theta-\theta\right)^2\right)
\end{eqnarray}
    where $c_n=\frac{1}{n}\sum_{i=1}^n\left[\dot{m}_{\theta}\,f_{\epsilon}\left(x-m_{\theta}\left(X_{i-1}\right)\right)\right]$.
By (\ref{thest}) and (\ref{new3}),
\begin{eqnarray*}
&&\frac{1}{\sqrt{n}}\left(\tilde{S}_n(x)-2nF_X(x)\right)\\
&=&\frac{1}{\sqrt{n}}\sum_{i=1}^n\left[F_{\epsilon}\left(x-m_{\theta}\left(X_{i-1}\right)\right)+F_g\left(x-\epsilon_i\right)-2F_X(x)+c_0R\left(X_{i-1},\theta,\epsilon_i\right)\right]\\
&&~~~~~~~~~~~~~+\left(c_n-c_0\right)\frac{1}{\sqrt{n}}\sum_{i=1}^nR\left(X_{i-1},\theta,\epsilon_i\right)-c_n\,o_{\mathbb{P}}(1)+\frac{1}{\sqrt{n}}\,O\left(\left(\hat\theta-\theta\right)^2\right)
\end{eqnarray*}
    Then, by $\frac{1}{\sqrt{n}}\sum_{i=1}^nR\left(X_{i-1},\theta,\epsilon_i\right)\,\rightarrow\,N\left(0,\|\sum_{i=1}^{\infty}\mathcal{P}_0R\left(X_{i-1},\theta,\epsilon_i\right)\|^2\right)$ and by (\ref{thest}), Lemma 5 and $c_n\,\,\stackrel{\mathbb{P}}{\rightarrow}\,\,c_0$, the lemma follows.
$$\qquad\qquad\qquad\qquad\qquad\qquad\qquad\qquad\qquad\qquad\qquad\qquad\qquad\qquad\qquad\qquad\qquad\Box$$

    We are now ready to state the main theoretic result:
		
{\bf Theorem 1A.} For each $x\in\mathcal{X}$,
\begin{eqnarray*}
\sqrt{n}\left(\hat{F}_c(x)-F_X(x)\right)\,\,\,\Rightarrow\,\,\,N\left(0,\,\,\sigma_2^2(x)\right)
\end{eqnarray*}
		
proof) Recall $\hat{F}_{c}(x)$ defined by (\ref{convol}), the convolution c.d.f. estimate when the innovation c.d.f. is unknown. By (\ref{eq38}), we have the following:
\begin{eqnarray}\label{eq39}
&&\int_{\mathbb{R}}\left[\hat{F}_{\epsilon}\left(x-y,\hat{\theta}\right)-F_{\epsilon}(x-y)\right]\,\left[\hat{f}_{g}\left(y,\hat{\theta}\right)-f_{g}(y)\right]dy\nonumber\\
&=&\hat{F}_{c}(x)-\int_{\mathbb{R}}F_{\epsilon}(x-y)\hat{f}_{g}\left(y,\hat{\theta}\right)dy-\int_{\mathbb{R}}\hat{F}_{\epsilon}\left(x-y,\hat{\theta}\right)f_g(y)dy+F_X(x)\nonumber\\
&=&O_{\mathbb{P}}\left(\sqrt{\frac{1}{n^2b^4}+\frac{b}{n}}\right)
\end{eqnarray}
    Note also
\begin{eqnarray}\label{new1}
\int_{\mathbb{R}}F_{\epsilon}(x-y)\hat{f}_{g}\left(y,\hat{\theta}\right)dy&=&\int_{\mathbb{R}}\frac{1}{nb}\sum_{i=1}^nK\left(\frac{y-m_{\hat\theta}\left(X_{i-1}\right)}{b}\right)F_{\epsilon}(x-y)dy
\end{eqnarray}
\begin{eqnarray}\label{new2}
\int_{\mathbb{R}}\hat{F}_{\epsilon}\left(x-y,\hat{\theta}\right)f_g(y)dy&=&\int_{\mathbb{R}}\frac{1}{n}\sum_{i=1}^nG\left(\frac{y-\hat\epsilon_i}{b}\right)f_g(x-y)dy
\end{eqnarray}
    Then, by (\ref{eq39})--(\ref{new2}),
\begin{eqnarray}\label{eq40}
\hat{F}_c(x)+F_X(x)&=&O_{\mathbb{P}}\left(\sqrt{\frac{1}{n^2b^4}+\frac{b}{n}}\right)+\frac{1}{n}\sum_{i=1}^n\int_{\mathbb{R}}K(u)F_{\epsilon}\left(x-ub-m_{\hat\theta}\left(X_{i-1}\right)\right)du\nonumber\\
&&~~~~~~~~~~~~~~~~~~~~~~~~~~+\frac{b}{n}\sum_{i=1}^n\int_{\mathbb{R}}G(u)f_g\left(x-ub-\hat\epsilon_i\right)du
\end{eqnarray}
    By the integration-by-parts,
\begin{eqnarray}\label{eq41}
\int_{\mathbb{R}}G(u)f_g\left(x-ub-\hat\epsilon_i\right)du&=&-\frac{1}{b}G(u)F_g\left(x-ub-\hat\epsilon_i\right)\Big|^{u=\infty}_{u=-\infty}+\frac{1}{b}\int_{\mathbb{R}}G^{\prime}(u)F_g\left(x-ub-\hat\epsilon_i\right)du\nonumber\\
&=&\frac{1}{b}\int_{\mathbb{R}}K(u)F_g\left(x-ub-\hat\epsilon_i\right)du
\end{eqnarray}
    By (\ref{eq40}) and (\ref{eq41}),
\begin{eqnarray}\label{eq42}
\hat{F}_c(x)+F_X(x)&=&O_{\mathbb{P}}\left(\sqrt{\frac{1}{n^2b^4}+\frac{b}{n}}\right)+\frac{1}{n}\int_{\mathbb{R}}K(u)\tilde{S}_n(x-ub)du
\end{eqnarray}
    By a Taylor's expansion,
\begin{eqnarray}\label{eq43}
\int_{\mathbb{R}}K(u)\tilde{S}_n(x-ub)du=\tilde{S}_n(x)+\frac{b^2\phi_K}{2}\tilde{S}^{\prime\prime}_n(x)+\frac{b^4}{4!}\int_{\mathbb{R}}u^4K(u)\tilde{S}^{(4)}_n(x_0)du
\end{eqnarray}
    where $\phi_K=\int_{\mathbb{R}}u^2K(u)du$ and $x_0\in(x,\,x-ub)$. Then, by (\ref{eq42}) and (\ref{eq43}),
\begin{eqnarray}\label{eq44}
&&\sqrt{n}\left(\hat{F}_c(x)-F_X(x)\right)\nonumber\\
&=&\frac{1}{\sqrt{n}}\left(\tilde{S}_n(x)-2nF_X(x)\right)+O_{\mathbb{P}}\left(\sqrt{\frac{1}{nb^4}+b}\right)\nonumber\\
&&~~~~~~~~~~~~~~~~~~~~~+\frac{\phi_Kb^2\sqrt{n}}{2}\left(\frac{\tilde{S}^{\prime\prime}_n(x)}{n}\right)+\frac{b^4\sqrt{n}}{4!}\left(\frac{1}{n}\int_{\mathbb{R}}u^4K(u)\tilde{S}^{(4)}_n(x_0)du\right)
\end{eqnarray}
    where $\tilde{S}^{\prime\prime}_n(x)=\sum_{i=1}^n\left[f^{\prime}_{\epsilon}\left(x-m_{\hat\theta}\left(X_{i-1}\right)\right)+f^{\prime}_g\left(x-\hat\epsilon_i\right)\right]$. Hence, by Assumption 3 and Lemma 6, the theorem follows.
$$\qquad\qquad\qquad\qquad\qquad\qquad\qquad\qquad\qquad\qquad\qquad\qquad\qquad\qquad\qquad\qquad\qquad\Box$$

{\bf Theorem 1B.} Let $D_2(x)=\sum_{i=1}^{\infty}\mathcal{P}_0\left[F_{\epsilon}\left(x-m_{\theta}(X_{i-1})\right)+F_g(x-\epsilon_i)+c_0R\left(X_{i-1},\theta,\epsilon_i\right)\right]$. Then, for any compact interval $\mathcal{X}$,
\begin{eqnarray*}
\left\{\sqrt{n}\left[\hat{F}_c(x)-F_X(x)\right],\,\,x\in\mathcal{X}\right\}\,\,\,\Rightarrow\,\,\,\left\{W_2(x),\,\,x\in\mathcal{X}\right\}
\end{eqnarray*}
  where $W_2(x)$ is a mean-zero Gaussian process with covariance function $cov\left[W_2(x),\,W_2(x')\right]=\mathbb{E}\left[D_2(x)D_2(x')\right]$.

proof) It is straightforward to verify the finite-dimensional convergence based on the Cram\'{e}r-Wold device. Hence we need to verify the tightness condition. By the Lipschitz continuity of $F_X(\cdot)$ and $F_{\epsilon}(\cdot)$ ({\it i.e.} Assumption 3),
\begin{eqnarray}\label{t0927}
&&\left\|\sum_{i=1}^n\left[F_{\epsilon}(x-m_{\hat\theta}(X_{i-1}))-F_X(x)\right]-\sum_{i=1}^n\left[F_{\epsilon}(x'-m_{\hat\theta}(X_{i-1}))-F_X(x')\right]\right\|^2\nonumber\\
&\leq&\sum_{i=1}^n\left\|\left[F_{\epsilon}(x-m_{\hat\theta}(X_{i-1}))-F_{\epsilon}(x'-m_{\hat\theta}(X_{i-1}))\right]\right\|^2+
\sum_{i=1}^n\left\|\left(F_X(x')-F_X(x)\right)\right\|^2\nonumber\\
&\leq&|x-x'|^2\,O(n)
\end{eqnarray}
Similarly, by the Lipschitz continuity of $F_X(\cdot)$ and $F_{g}(\cdot)$ ({\it i.e.} Assumption 3),
\begin{eqnarray}\label{t09272}
&&\left\|\sum_{i=1}^n\left[F_{g}(x-\hat{\epsilon}_i)-F_X(x)\right]-\sum_{i=1}^n\left[F_{g}(x'-\hat{\epsilon}_i)-F_X(x')\right]\right\|^2\nonumber\\
&\leq&\sum_{i=1}^n\left\|\left[F_{g}(x-\hat{\epsilon}_i)-F_{g}(x'-\hat{\epsilon}_i)\right]\right\|^2+\sum_{i=1}^n\left\|\left(F_X(x')-F_X(x)\right)\right\|^2\nonumber\\
&\leq&|x-x'|^2\,O(n)
\end{eqnarray}
By (\ref{eq44}),
\begin{eqnarray}\label{feb7}
&&\sqrt{n}\left(\hat{F}_c(x)-F_X(x)\right)\nonumber\\
&=&\frac{1}{\sqrt{n}}\sum_{i=1}^n\left[F_{\epsilon}(x-m_{\hat\theta}(X_{i-1}))-F_X(x)+F_{g}(x-\hat{\epsilon}_i)-F_X(x)\right]+o_{\mathbb{P}}(1)
\end{eqnarray}
Hence, by applying (\ref{t0927}) and (\ref{t09272}) to (\ref{feb7}),
\begin{eqnarray}\label{feb8}
\frac{\mathbb{E}\left|\sqrt{n}\left(\hat{F}_c(x)-F_X(x)\right)-\sqrt{n}\left(\hat{F}_c(x')-F_X(x')\right)\right|^2}{|x-x'|^2}&=&O(1)
\end{eqnarray}
Therefore, the tightness easily follows.
$$\qquad\qquad\qquad\qquad\qquad\qquad\qquad\qquad\qquad\qquad\qquad\qquad\qquad\qquad\qquad\qquad\qquad\Box$$

\section{Specification Test}

Given the kernel and convolution c.d.f. estimates in (\ref{kern}) and (\ref{convol}), we define:
\begin{eqnarray}\label{testat}
V_{T,b}(u):=\sqrt{n}\left(\hat{F}_{k}(u)-\hat{F}_{c}(u)\right)
\end{eqnarray}
where $\hat{F}_{k}(u)$ and $\hat{F}_{c}(u)$ are the kernel estimate in (\ref{kern}) and the convolution estimate in (\ref{convol}), respectively. We propose the following statistic:

{\bf Theorem 2.}
\begin{eqnarray}\label{tdist}
\int_{U}V^2_{T,b}(u)\pi(u)du\,\,\Rightarrow\,\,\int_{U}Z^2(u)\pi(u)du
\end{eqnarray}
where $\int_{U}\pi(u)du=1$ and $Z(\cdot)$ is a Gaussian process with covariance kernel $\sigma_v^2(u,\,u')$.

proof) (i) Note that
\begin{eqnarray*}
V_{T,b}(u)=\sqrt{n}\left(\hat{F}_{k}(u)-F_{X}(u)\right)-\sqrt{n}\left(\hat{F}_{c}(u)-F_{X}(u)\right)
\end{eqnarray*}
By the Liapunov C.L.T., for any fixed $u\in\,U$,
\begin{eqnarray}\label{knorm}
\sqrt{n}\left(\hat{F}_{k}(u)-F_{X}(u)\right)\,\,\Rightarrow\,\,N\left(0,\,\sigma_k^2(u)\right)
\end{eqnarray}
where $\sigma_k^2(u)=F_X(u)\left(1-F_X(u)\right)$. Moreover, by Theorem 1A,
\begin{eqnarray}\label{cnorm}
\sqrt{n}\left(\hat{F}_{c}(u)-F_{X}(u)\right)\,\,\Rightarrow\,\,N\left(0,\,\sigma_2^2(u)\right)
\end{eqnarray}
where $\sigma_2(u)=\left\|\sum_{i=1}^{\infty}\mathcal{P}_0\left[F_{\epsilon}\left(u-m_{\theta}(X_{i-1})\right)+F_g\left(u-\epsilon_i\right)+c_0R\left(X_{i-1},\theta,\epsilon_i\right)\right]\right\|$.
Hence, for any fixed $u\in\,U$,
\begin{eqnarray}\label{fre}
V_{T,b}(u)\,\,\Rightarrow\,\,N\left(0,\,\sigma_v^2(u,\,u)\right)
\end{eqnarray}
where $\sigma_v^2(u,\,u):=\sigma_k^2(u)+\sigma_2^2(u)-2\,C(u,u)$. Here the covariance kernel $C(u,u')$ is given by the limit of:
\begin{eqnarray*}
-\frac{1}{n^2b}\sum_{i=1}^n\sum_{j=1}^n\sum_{k=1}^n Cov\left(G\left(\frac{u-X_i}{b}\right),\,\int G\left(\frac{u'-t-X_j+m_{\hat\theta}(X_{j-1})}{b}\right)K\left(\frac{u'-m_{\hat\theta}(X_{k-1})}{b}\right)dt\right)
\end{eqnarray*}
as $n\rightarrow\infty$. The Cramer-Wold device ensures:
\begin{eqnarray}\label{cw}
\left(\begin{array}{c}
V_{T,b}(u)\\
V_{T,b}(u')
 \end{array}\right)\,\,\,\Rightarrow\,\,\,N\left(
\left(\begin{array}{c}
0\\
0
\end{array}\right),\,\,\,\,\,
\left(\begin{array}{cc}
  	\sigma_v^2(u,\,u) & \sigma_v^2(u,\,u') \\ \sigma_v^2(u,\,u') & \sigma_v^2(u',\,u')
  \end{array}\right)
\right)
\end{eqnarray}\\

(ii) It leaves us to prove the $tightness$ condition. Note that:
\begin{eqnarray*}
&&V_{T,b}(u)-V_{T,b}(u')\\
&=&\sqrt{n}\left(\hat{F}_{k}(u)-F_{X}(u)\right)-\sqrt{n}\left(\hat{F}_{k}(u')-F_{X}(u')\right)+\sqrt{n}\left(\hat{F}_{c}(u')-F_{X}(u')\right)-\sqrt{n}\left(\hat{F}_{c}(u)-F_{X}(u)\right)
\end{eqnarray*}
Hence
\begin{eqnarray}\label{t1}
\frac{\mathbb{E}\left|V_{T,b}(u)-V_{T,b}(u')\right|^2}{|u-u'|^2}&\leq&\frac{\mathbb{E}\left|\sqrt{n}\left(\hat{F}_k(u)-F_X(u)\right)-\sqrt{n}\left(\hat{F}_k(u')-F_X(u')\right)\right|^2}{|u-u'|^2}\nonumber\\
&&+\frac{\mathbb{E}\left|\sqrt{n}\left(\hat{F}_c(u')-F_X(u')\right)-\sqrt{n}\left(\hat{F}_c(u)-F_X(u)\right)\right|^2}{|u-u'|^2}
\end{eqnarray}
Note that
\begin{eqnarray}\label{t2}
\frac{\mathbb{E}\left|\sqrt{n}\left(\hat{F}_k(u)-F_X(u)\right)-\sqrt{n}\left(\hat{F}_k(u')-F_X(u')\right)\right|^2}{|u-u'|^2}=O(1)
\end{eqnarray}
Moreover, by (\ref{feb8}),
\begin{eqnarray}\label{t3}
\frac{\mathbb{E}\left|\sqrt{n}\left(\hat{F}_c(u)-F_X(u)\right)-\sqrt{n}\left(\hat{F}_c(u')-F_X(u')\right)\right|^2}{|u-u'|^2}&=&O(1)
\end{eqnarray}
By applying (\ref{t2}) and (\ref{t3}) to (\ref{t1}),
\begin{eqnarray}
\frac{\mathbb{E}\left|V_{T,b}(u)-V_{T,b}(u')\right|^2}{|u-u'|^2}\,=\,O(1)
\end{eqnarray}
Hence the tightness easily follows.

$$\qquad\qquad\qquad\qquad\qquad\qquad\qquad\qquad\qquad\qquad\qquad\qquad\qquad\qquad\qquad\qquad\qquad\Box$$

{\bf Remark 1.} Similarly, we can formulate a test statistic based on density estimates:
\begin{eqnarray}\label{testat2}
v_{T,b}(u):=\sqrt{nb}\left(\hat{f}_{k}(u)-\hat{f}_{c}(u)\right)
\end{eqnarray}
where $\hat{f}_{k}(u)$ is a kernel density estimate and $\hat{f}_{c}(u)$ is the convolution density estimate from Kim et al (2015). From the root-n convergence of the convolution density estimate,
$$v_{T,b}(u)\,\,\Rightarrow\,\,N\left(0,~~f(x)\int_{\mathbb{R}}K^2(u)du\right)$$
where $f(x)$ is the true density. Hence one can propose the following similar test statistics:
$$\int_{U}v^2_{T,b}(u)\pi(u)du\,\,\Rightarrow\,\,\int_{U}Z^2(u)\pi(u)du$$
where $Z(u)$ is a Gaussian process. Note, however, that the convergence rate in (\ref{testat}) is faster than that from (\ref{testat2}) given the order of the bandwidth $b$. Thus, our test based on (\ref{testat}) is expected to perform better than the test based on (\ref{testat2}).

\section{Simulation study}
\subsection{Setup}
In this section we consider testing for model specification described in Section 2. Consider hypothesis testing
\begin{equation}
H_{0}: X_{i}=\theta_0|X_{i-1}|+\epsilon_{i},\quad\quad H_{a}: X_{i}=\theta_0|X_{i-1}|+\epsilon_{i}\sqrt{ \theta_{0}^{2}+\theta_1^{2}X_{i-1}^{2}  }
\end{equation}
where $\epsilon\sim N(0,1)$. We generate random sample of $\{X_{i}:i=1,2,..., T\}$. When we generate the sample, we set  $\theta_0=0.3$, $\theta_1=0.9$ and $T=200,400,$ and $600$. We demonstrate that our proposed test outperforms the benchmark test-see Kim et al (2015). To that end, we report empirical levels and powers and compare the findings with those of the benchmark test. The benchmark test, however, did not report powers, and hence, we compute them by monte carlo simulation as described therein. For our proposed test we employ block-wise bootstrap method proposed by K$\ddot{\textrm{u}}$nsch (1989) and Liu and Singh (1992). We first determine $l_{B}$, a size of the block so that the number of blocks, $n_{B}$, is $T/l_{B}$. Naik-Nimbalakar and Rajarshi (1994) showed that weak convergence of block-wise bootstrapped empirical process depends on the order of the $l_{B}$. They obtained desired results when $l_{B}=O(n^{1/2}-\epsilon)$, with $0<\epsilon<\frac{1}{2}$. Motivated by their work, $l_{B}=\{10,15,20,25\}$ are tried; we found that the proposed test displays the optimal result when $l_{B}=10$ for all $T$. Once we determine the value of $l_{B}$, we construct a block: we draw any uniform random number between 1 and $T-l_{B}+1$, say $k$, and choose $l_{B}$ consecutive observations, $X_{k+1},...,X_{k+l_{B}}$. We repeat constructing a block $n_{B}$ times, combine these $n_{B}$ blocks all together, and obtain resampled observations, $X_{1}^{*},...,X_{T}^{*}$. Recall $V_{T,b}(u)$ in (\ref{testat}). Let $V_{T,b}^{2}$ denote the integral of $V_{T,b}^{2}(u)$ as in (\ref{tdist}). For the calculation of the statistics, we use uniform kernel function: $K(u):=2^{-1}I(|u|\leq 1)$ where $I(\cdot)$ is an indicator function. Therefore,
\begin{displaymath}
G(u)=\int_{-\infty}^{u}K(x)dx =\left\{
   \begin{array}{ll}
   0, & \hbox{$u<-1$;} \\
   \frac{u+1}{2}, & \hbox{$-1\leq u <1$;} \\
   1, & \hbox{$u\geq 1$.}
  \end{array}
\right.
\end{displaymath}
For $\hat{\theta}$, we use least squares estimator. Define $h(t) := (-t^2+2c_{i}t)/8b$ where $c_{i}:=x+b-\hat{\epsilon}_i$. Then $\hat{F}_{c}$ in (\ref{convol}) can be rewritten as
\begin{eqnarray*}
\hat{F}_{c}(x) &=& \frac{1}{n^{2}b}\sum_{i=1}^{n}\sum_{j=1}^{n}\int G\left(\frac{x-t-\hat{\epsilon}_i}{b}\right)K\left(\frac{x-m_{\hat{\theta}}\left(X_{j-1}\right)}{b}\right) \,dt\\
&=& \frac{1}{n^{2}b}\sum_{i=1}^{n}\sum_{j=1}^{n} IF_{ij}(x),\\
\end{eqnarray*}
where
\begin{eqnarray*}
IF_{ij}(x)=\left\{
             \begin{array}{ll}
               0, & \hbox{$x<m_{\hat{\theta}}\left(X_{j-1}\right)+ \hat{\epsilon}_i-2b $;} \\
               h(x-\hat{\epsilon}_i+b)-h(m_{\hat{\theta}}\left(X_{j-1}\right)-b) ,
                       & \hbox{$m_{\hat{\theta}}\left(X_{j-1}\right)+ \hat{\epsilon}_i-2b\leq x<m_{\hat{\theta}}\left(X_{j-1}\right)+ \hat{\epsilon}_i$;} \\
               bG\left(\frac{x-m_{\hat{\theta}}\left(X_{j-1}\right)-\hat{\epsilon}_i-b}{b}\right)+ h(m_{\hat{\theta}}\left(X_{j-1}\right)+b)&\hbox{}\\
           \quad\quad\quad\quad\quad\quad\quad - h(x-\hat{\epsilon}_i-b),    & \hbox{$m_{\hat{\theta}}\left(X_{j-1}\right)+ \hat{\epsilon}_i\leq x<m_{\hat{\theta}}\left(X_{j-1}\right)+ \hat{\epsilon}_i+2b$;} \\
               b, & \hbox{otherwise.}
             \end{array}
           \right.
\end{eqnarray*}
Consequently, the great deal of simplification of $V_{T,b}(u)$ in (\ref{testat}) follows directly.

Define the bootstrap test statistic
\begin{equation}
V_{T,b}^{2*} = \int_{U} \big( V_{T,b}^{*}(u)-V_{T,b}(u)  \big)^{2}\pi(u)\,du
\end{equation}
where $V_{T,b}^{*}(u)$ denotes the counterpart of $V_{T,b}(u)$ which is obtained from resampled observations. We repeat block-wise bootstrap $B_{Iter}$ times, obtain $V_{T,b}^{2*}$'s, and calculate $100(1-\alpha)$ percentiles, $q_{1-\alpha}^{*}$. As various $l_{B}$'s are tried, so are $B_{Iter}$'s. Our findings show that empirical levels approaches more closely to suggested significance level $\alpha$ as $B_{Iter}$ increases. See, e.g., Table \ref{table:Blvary}. After $q_{1-\alpha}^{*}$ is obtained, we reject $H_0$ if $V_{T,b}^{2} > q_{1-\alpha}^{*}$. As a final step, we repeat this procedure 1000 times, count the number of rejections, and obtain empirical levels and powers by dividing it by 1000.

\subsection{Selection of $B_{Iter}$, $b$, and, $l_{B}$}
In the simulation study, $b=\{0.05, 0.1, 0.15, 0.2\}$ are tried for bandwidth. Since the choice of $b$ does not affect the powers and levels much, we only report the result corresponding to $b=0.1$. Table \ref{table:Blvary} reports empirical levels corresponding to various sizes of block and numbers of bootstrap iteration. As shown in the table, we obtain the optimal result at $(l_{B}, B_{Iter})=(10,200)$.
\begin{table}
\begin{center}
\begin{tabular}{ l c c  c c c||  c c  c c c  }
\hline
     & \multicolumn{5}{c}{$l_{B}=8$} & \multicolumn{5}{c}{$l_{B}=10$} \\
\cline{1-11}
  $\alpha$ & $B_{Iter}=40$  & 80 & 120  & 160 & 200 & $B_{Iter}=40$  & 80 & 120  & 160 & 200 \\
\hline
 0.1   & 0.149 & 0.124 & 0.119 & 0.114 & 0.109 & 0.125 & 0.108 & 0.104 & 0.107 & 0.101  \\
 0.075 & 0.107 & 0.098 & 0.086 & 0.084 & 0.081 & 0.102 & 0.088 & 0.084 & 0.078 & 0.078  \\
 0.05  & 0.079 & 0.072 & 0.063 & 0.059 & 0.055 & 0.074 & 0.064 & 0.053 & 0.047 & 0.045   \\
 0.025 & 0.056 & 0.046 & 0.030 & 0.033 & 0.027 & 0.044 & 0.035 & 0.028 & 0.028 & 0.026  \\
 0.01  & 0.050 & 0.025 & 0.016 & 0.014 & 0.016 & 0.029 & 0.013 & 0.012 & 0.014 & 0.011  \\
\hline
     & \multicolumn{5}{c}{$l_{B}=16$} & \multicolumn{5}{c}{$l_{B}=20$} \\
\cline{1-11}
 $\alpha$   & $B_{Iter}=40$  & 80 & 120  & 160 & 200 & $B_{Iter}=40$  & 80 & 120  & 160 & 200 \\
\hline
 0.1   &  0.128 & 0.119 & 0.127& 0.124& 0.120 & 0.158& 0.142& 0.134& 0.130& 0.127 \\
 0.075 &  0.091 & 0.087 & 0.085& 0.086& 0.089 & 0.131& 0.107& 0.104& 0.103& 0.099\\
 0.05  &  0.064 & 0.057 & 0.055& 0.054& 0.052 & 0.097& 0.076& 0.075& 0.067& 0.065\\
 0.025 &  0.041 & 0.032 & 0.029& 0.025& 0.026 & 0.065& 0.051& 0.039& 0.035& 0.031\\
 0.01  &  0.025 & 0.015 & 0.015& 0.012& 0.013 & 0.058& 0.029& 0.020& 0.019& 0.016\\
\hline

\end{tabular}
\end{center}
\vspace{-.4cm}
\caption{Levels when $B_{Iter}$ and $l_{B}$ vary with $T$ being fixed at 400.}\label{table:Blvary}
\end{table}
Table \ref{table:comparison} compares the proposed test with the benchmark test. It is hard to tell which test is superior in terms of the level. However, there is no room for argument in terms of the power: the proposed test dominates the benchmark test. When $T=200$, the differences of the powers between two tests are more than 0.3 for all $\alpha's$. When $T$ increase, the differences decrease: approximately 0.15 (0.1) for all $\alpha$'s when $T$ is 400 (600). However, benchmark test does not obtain the power of 0.9 for most of all $\alpha$'s even though $T$ reaches 600; when $\alpha=0.01$, power is still smaller than 0.8. On the contrary, our proposed test accomplishes the power more than 0.9 except a few cases: $\alpha=0.05,\,0.025,\,0.01$ with $T=200$ and $\alpha=0.01$ with $T=400$. Therefore, we conclude that the proposed test is much superior to the benchmark test.


\begin{table}
\begin{center}
\begin{tabular}{c  l c c  c c c  c }
\hline
   & $\alpha$ & \multicolumn{2}{c}{T=200} & \multicolumn{2}{c}{T=400} & \multicolumn{2}{c}{T=600}\\
\cline{3-8}
 &  & $V_{T,b}^{2}$  & $V_{T,b}^{\sup}$ & $V_{T,b}^{2}$  & $V_{T,b}^{\sup}$ & $V_{T,b}^{2}$ & $V_{T,b}^{\sup}$ \\
\hline
\multirow{5}{*}{Level} & 0.1    & 0.102 & 0.097 & 0.099 & 0.104 & 0.115 & 0.098 \\
                       & 0.075  & 0.073 & 0.074 & 0.078 & 0.077 & 0.078 & 0.073 \\
                       & 0.05   & 0.041 & 0.052 & 0.045 & 0.050 & 0.058 & 0.054 \\
                       & 0.025  & 0.024 & 0.022 & 0.026 & 0.028 & 0.024 & 0.029 \\
                       & 0.01   & 0.013 & 0.007 & 0.013 & 0.010 & 0.012 & 0.007 \\
\hline
\multirow{5}{*}{Power} & 0.1    & 0.940 & 0.604 & 0.974 & 0.822 & 0.984 & 0.900 \\
                       & 0.075  & 0.909 & 0.567 & 0.966 & 0.804 & 0.982 & 0.884 \\
                       & 0.05   & 0.860 & 0.519 & 0.951 & 0.760 & 0.976 & 0.870 \\
                       & 0.025  & 0.783 & 0.434 & 0.911 & 0.712 & 0.952 & 0.836 \\
                       & 0.01   & 0.687 & 0.332 & 0.860 & 0.636 & 0.930 & 0.785 \\
\hline
\end{tabular}
\end{center}
\vspace{-.4cm}
\caption{Proposed test vs the benchmark test when $T=200$, 400, and 600}\label{table:comparison}
\end{table}

\end{document}